\documentclass[11pt]{amsart}
\usepackage{tabularx,booktabs,tikz}
\usepackage{caption}
\usepackage{amsmath}
\usepackage{amsfonts}
\usepackage{amscd}
\usepackage{amsthm}
\usepackage{amssymb} 
\usepackage{latexsym}
\usepackage{euscript}
\usepackage{epsfig}
\usepackage{graphics}
\usepackage{array}
\usepackage{enumerate}
\usepackage{dsfont}
\usepackage{color}
\usepackage{wasysym}
\usepackage{hyperref}
\usepackage{pdfsync}
\usepackage{enumitem}
\usepackage{fix-cm}
\usepackage{anyfontsize} 
\usepackage{wasysym}

\usepackage{newtxtext,newtxmath}
\usepackage[english]{babel}  
\newenvironment{remarklist}{
    \begin{enumerate}[
        leftmargin=*,
        itemsep=2pt,
        parsep=2pt,
        topsep=4pt
    ]
}{
    \end{enumerate}
}

\usepackage[left=3.5cm, right=2.5cm, top=3cm, bottom=3cm]{geometry}

\newcommand{\Hmm}[1]{\leavevmode{\marginpar{\tiny%
$\hbox to 0mm{\hspace*{-0.5mm}$\leftarrow$\hss}%
\vcenter{\vrule depth 0.1mm height 0.1mm width \the\marginparwidth}%
\hbox to
0mm{\hss$\rightarrow$\hspace*{-0.5mm}}$\\\relax\raggedright #1}}}

\newcommand{\Capp}[2][]{\mathrm{Cap}^{#1}_{#2}}

\newtheorem{theorem}{Theorem}[section]
\newtheorem{lemma}[theorem]{Lemma}

\newtheorem{definition}[theorem]{Definition}

\newtheorem{remark}[theorem]{Remark}

\newtheorem{example}[theorem]{Example}

\allowdisplaybreaks

\begin{document}

\title[$p$-isocapacitary constants associated with $p$-Laplacian on graphs]{$p$-isocapacitary constants associated with $p$-Laplacian on graphs}

\author{Bobo Hua}
\address{Bobo Hua: School of Mathematical Sciences, LMNS, Fudan University, Shanghai\ 200433, China}

\email{bobohua@fudan.edu.cn} 

\author{Lili Wang}
\address{Lili Wang: School of Mathematics and Statistics, Key Laboratory of Analytical Mathematics and Applications (Ministry of Education), Fujian Key Laboratory of Analytical Mathematics and Applications (FJKLAMA), Fujian Normal University, 350117 Fuzhou, P.R. China.}

\email{liliwang@fjnu.edu.cn}

\begin{abstract}

In this paper, we introduce $p$-isocapacitary constants for the $p$-Laplacian on graphs and apply them to derive estimates for the first eigenvalues of the Dirichlet $p$-Laplacian, the Neumann $p$-Laplacian, and the $p$-Steklov problem.

\end{abstract}
\par
\maketitle

\bigskip

\section{introduction}
Bounding the first eigenvalue of the $p$-Laplacian by isoperimetric constants is a central problem in differential geometry and graph theory. The standard approach for $p>1$, the Cheeger inequality, gives a lower bound of order $h^p$ and an upper bound of order $h$, where $h$ is the Cheeger constant. The bounds are of different orders, and the gap widens as $h \to 0$ \cite{Mazya1960,mazya1964,Mazya1985,mazya2009,Keller-Mugnolo16}. This paper replaces the Cheeger constant with a $p$-isocapacitary constant, obtaining bounds of matching order for the Dirichlet, Neumann, and Steklov eigenvalues of the $p$-Laplacian on graphs.

The idea of using capacity to estimate eigenvalues originates in the continuous setting. On a Riemannian manifold $M$, the $p$-capacity ($p>1$) of a capacitor $(F,\Omega)$---a compact set $F$ inside an open set $\Omega$---measures the minimal energy required to achieve electrostatic separation, defined as in \cite{grigor-isoperimetric} 
\begin{equation}\label{p-cap-continuous}
\operatorname{Cap}_p(F, \Omega) = \inf\limits_{u} \left\{ \int_\Omega |\nabla u|^p  d\mu  \;\bigg|\; 
u|_F = 1,  
u \in \text{Lip}_c(\Omega)
 \right\},  
\end{equation} where $\text{Lip}_c(\Omega)$ is the space of Lipschitz functions compactly supported in $\Omega.$ It yields Sobolev inequalities from the isoperimetric constant \cite{Mazya1960, Mazya1985}, characterizes the existence of $p$-harmonic functions on Riemannian manifolds \cite{Holopainen92}, controls the long-time behavior of the heat kernel \cite{barlow2009heat}, and bounds the first eigenvalue of the $p$-Laplacian from below \cite{Mazya1962, Mazya1985}. For recent applications to Schr\"{o}dinger equations and complex Hessian 
equations, see \cite{Barletta2025, Wang-Zhou2024}.

In the discrete setting, however, the standard tool remains the Cheeger inequality. On a weighted graph, the $p$-Laplacian is the natural counterpart of the continuous operator. Upper bounds follow from the variational characterization. Lower bounds are the principal difficulty. The Cheeger inequality originated on manifolds. Dodziuk, Alon--Milman, and Chung \cite{Dodziuk1984difference,1985lambda1,Chung1997SpectralGT} adapted the isoperimetric constant $h$ to graphs. For the $p$-Laplacian, Keller and Mugnolo \cite{Keller-Mugnolo16} proved
\[
C_1 h^p \leq \lambda_{1,p} \leq C_2 h,
\]
where $C_1,C_2$ depend on $p$. Consequently, the gap between the two bounds widens as $h \to 0$, precisely when the graph is nearly disconnected and a tight estimate of $\lambda_{1,p}$ is most relevant.

The isocapacitary constant provides an alternative. Maz'ya \cite{Mazya1962, mazya1964,mazya2009} introduced this constant in the continuum and established sharp two-sided eigenvalue bounds for the Dirichlet Laplacian. Hua et al. \cite{Hua-Munch-Wang} carried the approach over to graphs for $p=2$, obtaining analogous bounds for the graph Laplacian and for Steklov eigenvalues. We extend their framework to the $p$-Laplacian ($p>1$) by introducing a discrete $p$-isocapacitary constant. We also study the Steklov eigenvalue problem for the $p$-Laplacian on graphs, which has not been considered for $p \neq 2$ in the discrete setting. A discrete $p$-coarea formula expresses the $p$-energy as an integral of capacities over level sets. This formula is independent of the approach in \cite{Hua-Munch-Wang} and yields upper and lower bounds of the same order for Dirichlet, Neumann, and Steklov boundary conditions, overcoming the order mismatch inherent in the Cheeger inequality.

We now fix the notation for weighted graphs used throughout the paper.

Let $G=(V,E,w,m)$ be an undirected simple graph with vertex set $V$, edge set $E$, vertex weight $m: V\to \mathbb{R}_+$, and edge weight $w: V\times V\to \mathbb{R}_{\ge 0}$, written $w_{xy}=w(x,y)$, with $w_{xy}>0$ if $\{x,y\}\in E$ and $w_{xy}=0$ otherwise. Two vertices $x,y$ are neighbors, written $x\sim y$, if $\{x,y\}\in E$. For a subset $A\subset V$, its volume is $m(A)=\sum\limits_{x\in A} m(x)$. The graph is connected if for any $x,y\in V$ there exist vertices $z_0,\dots,z_n\in V$ such that $x=z_0\sim z_1\sim\cdots\sim z_n=y$. For subsets $\Omega,\Omega'\subset V$, the set of edges between $\Omega$ and $\Omega'$ is \[E(\Omega,\Omega'):=\{\{x,y\}\in E \mid x\in\Omega,\; y\in\Omega'\}.\]
We assume throughout that $G$ is locally finite, i.e., every vertex has finitely many neighbors. We call the quadruple $G=(V,E,w,m)$ a weighted graph.

For a subset $\Omega\subset V$, the vertex boundary of $\Omega$ is
\[
\delta\Omega:=\{y\in V\setminus\Omega \mid \exists\,x\in\Omega \text{ with } x\sim y\}.
\]
Set $\overline{\Omega}:=\Omega\cup\delta\Omega$. We assume that the induced subgraph on $\overline{\Omega}$ is connected.

For a set $A$, we denote by $\mathbb{R}^A$ the space of real-valued functions on $A$.
The $p$-Laplacian on graphs was introduced in \cite{Yamasaki1977,Yamasaki1986ideal}. For $\Omega\subset V$, $p\in(1,\infty)$, and $f\in\mathbb{R}^{\overline{\Omega}}$,
\[
\Delta_p^G f(x):=\frac{1}{m(x)}\sum_{y\in\overline{\Omega}} w_{xy}\,|f(y)-f(x)|^{p-2}\bigl(f(y)-f(x)\bigr),\qquad x\in\Omega.
\]
When the graph $G$ is clear from the context, we simply write $\Delta_p f(x)$.

For $X\subset V$, we define the graph $G_X=(\overline{X},E(X,\overline{X}),w,m)$, where $w$ and $m$ are restricted to $E(X,\overline{X})$ and $\overline{X}$, respectively. Note that edges within $\delta X$ are excluded from the edge set of $G_X$. In particular, $w_{xy}=0$ for $\{x,y\}\in E(\delta X,\delta X)$.

Next, we will study the Dirichlet and Neumann eigenvalue problems for the $p$-Laplacian on graphs, as well as the $p$-Steklov eigenvalue problem. 

We recall a well-known result of Maz'ya in the continuous setting: the first Dirichlet eigenvalue of the p-Laplacian is estimated by the p-capacity.
\begin{theorem}[\cite{Mazya1960,mazya1964,Mazya1985,mazya2009}]
For a bounded domain $\Omega$ in a Riemannian manifold,
\[
c_p\alpha_p^D(\Omega)\leq \lambda_{1,p}(\Omega)\leq \alpha_p^D(\Omega),
\]
where $\lambda_{1,p}$ is the first eigenvalue of the Dirichlet $p$-Laplacian on $\Omega$,
\[
\alpha^D_p(\Omega)=\inf\limits_{F\subset\subset \Omega}\frac{\mathrm{Cap}_p(F,\Omega)}{\mathrm{vol}(F)},
\]
and $c_p=(p-1)^{p-1}p^{-p}$.
\end{theorem}

We now turn to the discrete setting. For $\Omega\subset V$, the Dirichlet problem on $\Omega$ is
\begin{align*}
\begin{cases}
\Delta_p f(x)=-\lambda |f(x)|^{p-2}f(x), & x\in \Omega,\\[2pt]
f(x)=0, & x\in\delta\Omega.
\end{cases}
\end{align*}
We denote by $\lambda_{1,p}(\Omega)$ the first eigenvalue of this problem \cite{Hua-Wang20}. We define the Dirichlet $p$-isocapacitary constant as
\begin{equation}\label{Dirichlet-isocap}
\alpha_{p}^D(\Omega):=\inf\limits_{A\subset \Omega} \frac{\Capp[\Omega]{p}(A,\delta\Omega)}{m(A)},
\end{equation}
where $\Capp[\Omega]{p}(A,\delta\Omega)$ is a discrete analog of $p$-capacity; see \eqref{p-capacity-def} for the definition.

Motivated by the Cheeger-type inequalities for the $2$-capacity on finite graphs established in \cite{Hua-Munch-Wang}, we establish analogous bounds for the $p$-capacity with $p>1$. Throughout, we assume $p\in(1,\infty)$ and define the constant
\begin{equation}\label{C(p)}
C_p=p\ln 4+\left(2-2^{\frac{1}{1-p}}\right)^{1-p}.
\end{equation}

We first bound the first eigenvalue of the Dirichlet $p$-Laplacian on finite graphs in terms of the Dirichlet $p$-isocapacitary constant. The following result is a discrete analog of Maz'ya's theorem.
\begin{theorem}\label{capacity-est-Dcase}
Let $G$ be a weighted graph and $\Omega\subset V$ a finite subset. Then
\begin{equation*}
\frac{1}{2^p C_p}\,\alpha_p^D(\Omega)\leq \lambda_{1,p}(\Omega)\leq \alpha_p^D(\Omega).
\end{equation*}
\end{theorem}

We now turn to infinite graphs. Following the exhaustion approach established in \cite{Hua-Munch-Wang}, let $\{W_i\}_{i=1}^\infty$ be an exhaustion of $G$ as in Definition~\ref{exhaustion}. The bottom of the spectrum of the $p$-Laplacian on $G$ is given by
\[
\lambda_{1,p}(G) = \lim_{i\to\infty} \lambda_{1,p}(W_i),
\]
where the limit exists and is independent of the choice of exhaustion (see Section~4 for a proof).
We define the Dirichlet $p$-isocapacitary constant of $G$ by
\[
\alpha_p^D(G)=\inf\limits_{\substack{A\subset V\\|A|<+\infty}}\frac{\mathrm{Cap}_p(A)}{m(A)},
\qquad
\mathrm{Cap}_p(A)=\lim\limits_{i\to\infty}\mathrm{Cap}_p^{W_i}(A).
\]

Combining the exhaustion approach with Theorem~\ref{capacity-est-Dcase}, we now establish two-sided estimates for $\lambda_{1,p}(G)$, the bottom of the spectrum of the $p$-Laplacian on infinite graphs.
\begin{theorem}\label{Dirichlet-eigenvalue-capacity-estimate-on-infinite-graph}
Let $G$ be an infinite weighted graph. Then
\begin{equation*}
\frac{1}{2^p C_p}\,\alpha_p^D(G)\leq \lambda_{1,p}(G)\leq \alpha_p^D(G).
\end{equation*}
\end{theorem}

\begin{remark}
For an infinite graph, our estimate shows that $\lambda_{1,p}(G)$ and $\alpha_p^D(G)$ are of the same order, which is better than Cheeger inequality for normalized Dirichlet $p$-Laplacian in \cite{Keller-Mugnolo16}:
\[
\frac{2^{p-1}}{p^p}\left(h^D(G)\right)^p\leq \lambda_{1,p}(G)\leq h^D(G),
\]
where
\[
h^D(\Omega)=\inf\limits_{W\subset \Omega}\frac{|\partial W|_w}{m(W)}, \ \ \ |\partial W|_w:=\sum\limits_{\{x,y\}\in E(W,W^c)}w_{xy}.
\]
\end{remark}

For a subset $W\subset V,$ we denote by $$\mathcal{P}(W):=\{\{A,B\}: A,B\subset W, A\neq \emptyset, B\neq \emptyset\}$$ the set of pairs of nonempty subsets of $W.$  

Cheeger type estimates are well-established for the first nonzero eigenvalue of Neumann $p$-Laplacian in both continuous and discrete settings. However, the corresponding isocapacitary inequalities have not been developed. In this section, we establish a Neumann $p$-isocapacitary inequality.

For a finite subset $\Omega\subset V$, the Neumann problem on $\Omega$ is
\begin{align}
\begin{cases}\label{Neumann-Laplacian-equation}
\Delta_p f(x)=-\mu |f(x)|^{p-2}f(x), & x\in\Omega,\\[2pt]
|\nabla f|^{p-2}\dfrac{\partial f}{\partial n}(x)=0, & x\in\delta\Omega,
\end{cases}
\end{align}
where the $p$-th outward normal derivative at $z\in\delta\Omega$ is defined by
\begin{equation}\label{Neumann-bdy-cond}
\left(|\nabla f|^{p-2}\frac{\partial f}{\partial n}\right)(z):=\frac{1}{m(z)}\sum_{x\in \Omega} w_{xz}\,|f(z)-f(x)|^{p-2}\bigl(f(z)-f(x)\bigr).
\end{equation}
We denote the first nonzero eigenvalue of \eqref{Neumann-Laplacian-equation} by $\mu_{1,p}(\Omega)$. When $\delta\Omega=\emptyset$, i.e., $\Omega=V$, the problem \eqref{Neumann-Laplacian-equation} reduces to the eigenvalue problem for the $p$-Laplacian on a graph without boundary.

Inspired by the case $p=2$ in \cite{Hua-Munch-Wang}, we define the Neumann $p$-isocapacitary constant as
\begin{equation}\label{Neumann-isocapacity}
\alpha_p^N(\Omega)=\inf_{\{A,B\}\in\mathcal{P}(\Omega)}\frac{\Capp[\Omega]{p}(A,B)}{m(A)\wedge m(B)},
\end{equation}
where $a\wedge b:=\min\{a,b\}$ and $\Capp[\Omega]{p}(A,B)$ is defined in \eqref{p-capacity-def}. 
 
We now state the Neumann analog of Theorem~\ref{capacity-est-Dcase}.

\begin{theorem}\label{Neumann-capacity-estimate}
Let $G$ be a weighted graph and $\Omega\subset V$ a finite subset with $|\Omega|\geq 2$. Then
\begin{equation*}
\frac{1}{2^p C_p}\,\alpha_p^N(\Omega)\leq \mu_{1,p}(\Omega)\leq 2^{p-1}\,\alpha_p^N(\Omega).
\end{equation*}
\end{theorem}
\begin{remark}
\begin{remarklist}
 
 \item In case that $\Omega=V,$ the estimate of above theorem yields the $p$-isocapacitary estimate for $p$-Laplacian of a graph without boundary. For a finite connected weighted graph $G=(V,E,w,m),$
    \begin{equation}\label{isocap-est-for-nonbdy}
\frac{1}{2^pC_p}\alpha_p(G)\leq \mu_{1,p}(G)\leq 2^{p-1}\alpha_p(G),\end{equation} where $C_p$ is given by \eqref{C(p)} and $\mu_{1,p}(G)$ is the first nonzero eigenvalue of the $p$-Laplacian, and $$\alpha_p(G):=\inf\limits_{A,B\in\mathcal{P}(V)}\frac{\Capp[V]{p}(A,B)}{m(A)\wedge m(B)}.$$

\item Hua-Huang \cite{Hua-Huang2018} estimated the first nonzero eigenvalue of Neumann Laplacian by Cheeger constant  for $p=2$, and for general $p$, Keller-Mugnolo proved the following Cheeger estimate for a finite graph $G$ without boundary:
    \[
\frac{2^{p-1}}{p^p}\left(h(G)\right)^p\leq \mu_{1,p}(G)\leq 2^{p-1}h(G),
\]
where $h(G)$ is the Cheeger constant of $G$.
Note that our estimate is better than the above results with matching orders for upper and lower bounds in terms of geometric quantities.
\end{remarklist}
\end{remark}

Existing research on discrete $p$-Steklov eigenvalues has so far been restricted to the case $p=2$ (see \cite{HHW2017,hassannezhad-miclo20}), and a general theory for $p>1$ remains open. Motivated by the continuous setting in \cite{Bonder2001, Garcia-Azorero-Manfredi-Peral-2006, Verma2020, provenzano2022upper}, we extend the Steklov problem to $p>1$ on graphs: for a finite subset $\Omega\subset V$, the $p$-Steklov problem on $\Omega$ is
\begin{align*}
\begin{cases}
\Delta_p f(x)=0, & x\in\Omega,\\[2pt]
|\nabla f|^{p-2}\dfrac{\partial f}{\partial n}=\sigma |f|^{p-2}f, & x\in\delta\Omega.
\end{cases}
\end{align*}
We denote the first nonzero eigenvalue by $\sigma_{1,p}(\Omega)$.

We now introduce the Steklov $p$-isocapacitary constant. For $|\delta\Omega|\geq 2$, define
\begin{equation}\label{def-alpha-S-isocap}
\alpha_p^S(\Omega)=\inf_{\{A,B\}\in\mathcal{P}(\delta\Omega)}\frac{\Capp[\Omega]{p}(A,B)}{m(A)\wedge m(B)}.
\end{equation}

The following result bounds $\sigma_{1,p}(\Omega)$ in terms of $\alpha_p^S(\Omega)$.
\begin{theorem}\label{Steklov-capacity-estimate}
Let $G$ be a weighted graph and $\Omega\subset V$ a finite subset with $|\delta\Omega|\geq 2$. Then
\begin{equation*}
\frac{1}{2^p C_p}\,\alpha_p^S(\Omega)\leq \sigma_{1,p}(\Omega)\leq 2^{p-1}\,\alpha_p^S(\Omega).
\end{equation*}
\end{theorem}

\begin{remark}
Open Problem~4.6 in the survey of Colbois et al.\ \cite{Colbois2024} asks whether, on a compact Riemannian manifold, there exists an isoperimetric-type constant that yields two-sided bounds of the same order for the first Steklov eigenvalue. On graphs, Hua, M\"{u}nch, and Wang \cite{Hua-Munch-Wang} answered this question affirmatively for $p=2$  using isocapacitary constants. Theorem~\ref{Steklov-capacity-estimate} answers the same question on graphs for all $p>1$.
\end{remark}

The article is organized as follows. Section~\ref{preliminaries} gives the necessary background on the $p$-Laplacian and isocapacitary constants. Section~\ref{coarea-a=2} develops a discrete coarea formula that connects $p$-energy functionals to level set capacities. This is the foundation for deriving the eigenvalue bounds. 
Section~\ref{Dirichlet case} establishes two-sided bounds for $\lambda_{1,p}$ in terms of the Dirichlet $p$-isocapacitary constant $\alpha_p^D$, for both finite graphs (Theorem~\ref{capacity-est-Dcase}) and infinite graphs (Theorem~\ref{Dirichlet-eigenvalue-capacity-estimate-on-infinite-graph}). 
Section~\ref{The first Neumann eigenvalues of the $p$-Laplacian} extends these results to the Neumann case using reweighted graph sequences. 
Section~\ref{First Steklov Eigenvalue of $p$-Laplacian} addresses the $p$-Steklov problem by harmonically extending boundary functions and adapting the capacity-based estimates.

\section{preliminaries}\label{preliminaries}

Let $G=(V,E,w,m)$ be a weighted graph and $S\subset V$ be a finite subset. For a function $f\in\mathbb{R}^S$, we define the $l^p$ norm of $f$ as
\[
\|f\|_{p,S}=\left(\sum\limits_{x\in S}|f(x)|^pm(x)\right)^\frac{1}{p}.
\]
We also define the $l^\infty$ norm of $f$ as
\[
\|f\|_{l^\infty,S}:=\sup\limits_{x\in S}|f(x)|.
\]
The space of $l^p$ summable functions on $S$ is given by 
\[
l^p(S):=\{f\in\mathbb{R}^S: \|f\|_{p,S}<+\infty\}.
\]

For any subset $X\subset V$, we define $l_0(X)$ as the set of functions on $X$ with finite support. 
Given functions $f,g\in\mathbb{R}^{\overline{X}}$, we define
\[
\langle f,g\rangle_{\overline{X}}:=\sum\limits_{x\in \overline{X}}f(x)g(x)m(x)
\]
and  
\begin{equation}
\mathcal{E}_p^X(f,g)
=\sum\limits_{\{x,y\}\in E(X,\overline{X})}w_{xy}|f(y)-f(x)|^{p-2}\left(f(y)-f(x)\right)\left(g(y)-g(x)\right),    
\end{equation}
whenever the summation absolutely converges. 
Furthermore, for any $A,B\subset \overline{X}$, we define the $p$-capacity as 
\begin{equation}\label{p-capacity-def}
 \Capp[X]{p}(A,B)=\inf\{\mathcal{E}_p^X(f,f): f|_A=1,\  f|_B=0,\  f\in l_0(\overline{X})\}.   
\end{equation}
Clearly, $\Capp[X]{p}(A,B)=0$ for $A=\emptyset$ or $B=\emptyset$, and $\Capp[X]{p}(A, B)$ is monotone increasing in both $A$ and $B$. 

Furthermore, for a finite set $X$, the infimum in the definition of $\mathrm{Cap}_p^X(A,B)$ is attained by a unique function $f\in l_0(\overline{X})$. Indeed, since $\overline{X}$ is finite, the admissible class $\{f\in l_0(\overline{X}): f|_A = 1,\ f|_B = 0\}$ is a nonempty closed affine subset of the finite-dimensional space $\mathbb{R}^{\overline{X}}$. The $p$-energy functional $\mathcal{E}_p^X(\cdot,\cdot)$ is strictly convex 
for $p>1$. Moreover, the boundary conditions together with the nonemptiness and finiteness of $A$ and $B$ make the functional coercive on this affine set. By the direct method of the calculus of variations in finite dimensions, 
existence and uniqueness of the minimizer follow. 
This minimizer is the unique solution to the system
\[
\Delta_p^{G_X} f = 0 \quad \text{for } x\in \overline{X}\setminus \{A\cup B\},
\quad \text{with } f|_A = 1,\ f|_B = 0.
\]

Exhaustion by finite subsets is a fundamental concept in the study of infinite graphs.
\begin{definition} \label{exhaustion}
Let $G=(V,E,w,m)$ be an infinitely weighted graph. A sequence of subsets of vertices $\mathcal{W}=\{W_i\}_{i=1}^\infty$ is called an exhaustion of $G$, written as $\{W_i\}\uparrow V$, if it satisfies
\begin{enumerate}
\item $W_1\subset W_2 \subset \cdots \subset W_i\subset \cdots \subset V$;
    \item $|W_i|<+\infty,$ for all $i=1,2,\cdots$;
    \item $V=\bigcup\limits_{i=1}^\infty W_i$.
\end{enumerate}
\end{definition}

For any infinite subset $U \subset V$ with closure $\overline{U} = U \cup \delta U$, the $p$-capacity of a finite set $A \subset \overline{U}$ is defined as
\[
\Capp[U]{p}(A)=\inf\limits_f\{\mathcal{E}_p(f,f)| f|_A=1, f\in l_0(\overline{U})\}.
\]
When $U = V$, we simplify the notation to
\[
\Capp[]{p}(A):=\Capp[V]{p}(A).
\]
Let $U \subset V$ be an infinite subset and let $A \subset \overline{U}$ be a finite set. For any exhaustion $\{W_k\} \uparrow \overline{U}$, we have
\[
\Capp[U]{p}(A)=\lim\limits_{i\to \infty}\Capp[U]{p}(A,\delta_UW_i).
\]
In particular, for the full graph ($U = V$) and any exhaustion $\{W_i\} \uparrow V$,
\[\Capp[]{p}(A)=\lim\limits_{i\to \infty}\Capp[W_i]{p}(A),\]
where $\text{Cap}_p^{W_k}(A)$ is the $p$-capacity of $A$ in the finite subgraph induced by $W_k$.

The analysis of $p$-capacities uses discrete analogs of classical identities; the following generalized Green’s formula is crucial.
\begin{lemma}[Green's Formula]\label{lemma:green-formula}
    Let $\Omega \subset V$ be a finite subset. For any functions $f, g \in \mathbb{R}^{\overline{\Omega}}$, the following identity holds:
    \begin{equation}
        -\langle \Delta_p f, g\rangle_{\Omega} + \left\langle |\nabla f|^{p-2}\frac{\partial f}{\partial n}, g \right\rangle_{\delta \Omega} = \mathcal{E}_p^{\Omega}(f, g).
    \end{equation}
\end{lemma}

\begin{proof}
  Retain the standing assumptions on the graph and the notation
  $\overline{\Omega}=\Omega\cup\delta\Omega$.
  For brevity set
  \[
    A_{xy}:=\left| f(y)-f(x)\right|^{p-2}\bigl(f(y)-f(x)\bigr)w_{xy},
  \]
  for $x,y\in\overline{\Omega}$. Clearly, $A_{yx}=-A_{xy}$.
The $p$-energy can be reshaped using the antisymmetry of $A_{xy}$ 
 
\begin{align*}
\mathcal{E}_p^\Omega(f,g)
&=\sum_{\{x,y\}\in E(\Omega,\overline{\Omega})}A_{xy}\bigl(g(y)-g(x)\bigr) \\
&=\sum_{\{x,y\}\in E(\overline{\Omega},\overline{\Omega})\setminus E(\delta\Omega,\delta\Omega) }(A_{xy}g(y)+A_{yx}g(x))\\
&=\frac{1}{2}\sum_{x,y\in\overline{\Omega}}\bigl(A_{xy}g(y)+A_{yx}g(x)\bigr)-\frac{1}{2}\sum_{x,y\in\delta\Omega}\bigl(A_{xy}g(y)+A_{yx}g(x)\bigr) \\
&=\sum_{x,y\in\overline{\Omega}}A_{xy}g(y)-\sum_{x,y\in\delta\Omega}A_{xy}g(y)\\
&=\sum_{y\in\Omega}g(y)\sum_{x\in\overline{\Omega}}A_{xy}+\sum_{y\in\delta\Omega}g(y)\sum\limits_{x\in \Omega}A_{xy} .
\end{align*}
By definition of $\Delta_p$ and \eqref{Neumann-bdy-cond} we have
\[\sum\limits_{x\in\overline{\Omega}}A_{xy}=-m(y)\Delta_p f(y)\]
for $y\in\Omega$ and
\[
\sum_{x\in\Omega}A_{xy}=|\nabla f|^{p-2}\frac{\partial f}{\partial n}(y)
\]
for $y\in\delta\Omega$.
Therefore,
\[
\mathcal{E}_p^\Omega(f,g)=\langle\Delta_p f,g\rangle_\Omega+\Bigl\langle |\nabla f|^{p-2}\frac{\partial f}{\partial n},\,g\Bigr\rangle_{\delta\Omega}.
\]
\end{proof}

\section{Coarea formula}\label{coarea-a=2}

In this section, we establish a discrete coarea formula that serves as a fundamental tool for our subsequent analysis. We need the following key analytic inequality.

\begin{lemma}\label{Optimal constant inequality}
Let $a > 1$ and $u, v \geq 0$ satisfy $u \geq a v$. Then 
\[
\frac{u^p}{a^p} - v^p \leq \frac{(u - v)^p}{\left(a^\frac{p}{p-1}  -1\right)^{p-1}}.
\]
The above estimate is sharp.
\end{lemma}

\begin{proof}
If $v = 0$, then $u\geq 0$ and the inequality holds (and becomes an equality only if $u=0$). For $u>0$, it is sufficient to prove 
\[
\left(a^\frac{p}{p-1}  -1\right)^{p-1}\leq a^p=\left(a^\frac{p}{p-1}\right)^{p-1}. 
\]
It is clear from $a>1$.  

For the case $v > 0$,  
let $t = \frac{u}{v}$. Then $t \geq a$. It is sufficient to prove 
\[
\frac{t^p - a^p}{a^p (t - 1)^p} \leq \frac{1}{\left(a^{\frac{p}{p-1}} - 1\right)^{p-1}}.  
\]
Define 
\[
f(t) := \frac{t^p - a^p}{a^p (t - 1)^p}, \quad t \geq a.
\]
Thus,
\[
f'(t) = \frac{p\left(a^{p} - t^{p-1}\right)}{a^{p}(t-1)^{p+1}}.
\]
Analyzing $f'(t)$ indicates that $f(t)$ attains its maximum at $t = a^{\frac{p}{p-1}}$ and
\[
f(a^{\frac{p}{p-1}})=\frac{1}{\left(a^{\frac{p}{p-1}} - 1\right)^{p-1}}.
\]
This completes the proof.
\end{proof}

We now prove a discrete coarea formula, motivated by its continuous analog in \cite{Mazya2005Conductor}.
\begin{lemma}\label{coarea formular} 
Let $f\in \mathbb{R}^{\overline{\Omega}}$, $M_t=\{x\in \overline{\Omega}: |f(x)|\geq t\}$. For any constant $a>1$ and $1< p<\infty$, we have
\begin{equation}
\int_0^\infty \Capp[\Omega]{p}(M_{at},M_t^c)pt^{p-1}dt\leq C(a,p)\mathcal{E}^\Omega_p(f,f),
\end{equation}
where 
\[
C(a,p)=\frac{2p\ln a}{(a-1)^p}+2\left(a^\frac{p}{p-1}  -1\right)^{1-p}.
\]

In particular, 
\begin{equation*}
\int_0^\infty \Capp[\Omega]{p}(M_{2t},M_t^c)pt^{p-1}dt\leq C_p\mathcal{E}_p^{\Omega}(f,f),
\end{equation*}
where $C_p$ is given by \eqref{C(p)}.   
\end{lemma}

\begin{proof}
Let 
\begin{equation*}
\phi(x)=\frac{ (|f|-t)_+\wedge (a-1)t }{(a-1)t}. 
\end{equation*} 
Then
\begin{align*}\begin{cases}
\phi(x)\equiv 1, \ \ & x\in M_{at};\\
0\leq \phi(x)\leq 1, &x\in \Omega;\\
\phi(x)\equiv 0, & x\in M_t^c.
\end{cases}\end{align*}
By the definition of the capacity,  
\begin{equation*}
\Capp[\Omega]{p}(M_{at},M_t^c)\leq \mathcal{E}^\Omega_p(\phi,\phi).   
\end{equation*}
Thus,
\begin{equation}\label{key-inequality-1}
t^p\Capp[\Omega]{p}(M_{at},M_t^c)\leq \mathcal{E}^\Omega_p(t\phi,t\phi).
\end{equation}
Let $K_t=M_t\setminus M_{at}$. Since 
\[
\mathcal{E}^\Omega_p(t\phi,t\phi)=\sum\limits_{\{x,y\}\in E(\Omega,\overline{\Omega})}w_{xy}|t\phi(y)-t\phi(x)|^p
\]
and
\begin{align*}
\left|t\phi(y)-t\phi(x)\right|=
\begin{cases}
\frac{\big||f(y)|-|f(x)|\big|}{a-1}, \ \ \ \ \ \ \ \ \ \ \ \  &(x,y)\in K_t\times K_t;\\
\frac{at-|f(x)|}{a-1}\leq \frac{|f(y)|-|f(x)|}{a-1}, &
(x,y)\in  K_t\times M_{at};\\
\frac{|f(x)|-t}{a-1}\leq \frac{|f(x)|-|f(y)|}{a-1}, &
(x,y)\in K_t\times M^c_t;\\
t, &\{x,y\}\in E(M_{at}, M^c_t),
\end{cases}
\end{align*}
we have
\begin{align*} 
\mathcal{E}_p^\Omega(t\phi,t\phi)
&=\sum\limits_{\{x,y\}\in E(\Omega,\overline{\Omega})}w_{xy}|t\phi(x)-t\phi(y)|^p\\
&= \left(\sum\limits_{\{x,y\}\in E(K_t, \overline{\Omega})}   
 +\sum\limits_{\{x,y\}\in E(M_{at}, M^c_t)}\right) \left(w_{xy} |t\phi(y)-t\phi(x)|^p\right)\\
 &\leq\frac{1}{(a-1)^p}\sum\limits_{\{x,y\}\in E(K_t, \overline{\Omega})}w_{xy}|f(x)-f(y)|^p+\sum\limits_{\{x,y\}\in E(M_{at}, M^c_t)}w_{xy}t^p.
\end{align*}
Hence, \eqref{key-inequality-1} implies that
\begin{equation}
\begin{aligned}\label{key-inequality-12}
\int_0^\infty\Capp[\Omega]{p}(M_{at},M_t^c)pt^{p-1}dt&=p\int_0^\infty \frac{1}{t}\mathcal{E}^\Omega_p(t\phi,t\phi)dt\\
&\leq\frac{p}{(a-1)^p}\int_0^\infty\frac{1}{t}\sum\limits_{\{x,y\}\in E(K_t, \overline{\Omega})}w_{xy}|f(x)-f(y)|^p dt\\
&\ \ \ \ \ \ \ \ +\int_0^\infty \sum\limits_{\{x,y\}\in E(M_{at}, M^c_t)}w_{xy}d(t^p)\\
&:=\frac{p}{(a-1)^p}\text{I}+\text{II}.
\end{aligned}
\end{equation}

Since the integrand is nonnegative, Tonelli's theorem \cite[pages 67-68, Theorem 2.37 (a)]{Folland99} allows us to interchange the sum and the integral. Using this, we can estimate the first term as follows:
\begin{equation}\begin{aligned}\label{int-1}
\text{I}&\leq \int_0^\infty\frac{1}{t}\sum_{(x,y)\in K_t\times\overline{\Omega}}w_{xy}|f(x)-f(y)|^p dt\\
&=\int_0^\infty \frac{1}{t}\sum_{x,y\in \overline{\Omega} }w_{xy}\chi_{K_t}(x)|f(x)-f(y)|^p dt\\
&=\sum_{x,y\in \overline{\Omega} }w_{xy}|f(x)-f(y)|^p \int_0^\infty \frac{1}{t}\chi_{K_t}(x) dt\\
&=\sum_{x,y\in \overline{\Omega} }w_{xy}|f(x)-f(y)|^p \int_\frac{|f(x)|}{a}^{|f(x)|}\frac{1}{t}dt\\
&=2\ln a \sum_{{x,y}\in E(\Omega, \overline{\Omega})}w_{xy}|f(x)-f(y)|^p.
\end{aligned}
\end{equation}
By Lemma \ref{Optimal constant inequality} and $|f(x)|\geq a|f(y)|$, we have 
\[
\frac{|f(x)|^p}{a^p}-|f(y)|^p
\leq \left(a^\frac{p}{p-1}  -1\right)^{1-p} \left(|f(x)|-|f(y)|\right)^p\leq \left(a^\frac{p}{p-1}  -1\right)^{1-p}|f(x)-f(y)|^p. 
\]
Using Tonelli's theorem to interchange the sum and the integral, we can estimate $\operatorname{II}$ as follows:

\begin{align}\begin{split}\label{int-2}
\operatorname{II}=\int_0^\infty\sum\limits_{(x,y)\in M_{at}\times M_t^c}w_{xy}pt^{p-1}dt
&=\sum\limits_{(x,y)\in \overline{\Omega}\times\overline{\Omega}}w_{xy}\int_0^\infty\mathcal{X}_{M_{at}}(x)\mathcal{X}_{M_t^c}(y) pt^{p-1}dt\\
 &=\sum\limits_{x,y\in \overline{\Omega}}w_{xy}\int_{|f(y)|}^\frac{|f(x)|}{a}pt^{p-1}dt\\
 &\leq \left(a^\frac{p}{p-1}-1\right)^{1-p}\sum\limits_{x,y\in \overline{\Omega}}w_{xy}|f(x)-f(y)|^p\\
 &=2\left(a^\frac{p}{p-1}-1\right)^{1-p}\sum\limits_{\{x,y\}\in E(\Omega,\overline{\Omega})}w_{xy}|f(x)-f(y)|^p,
\end{split}\end{align} 
where the inequality follows from evaluating the integral and applying the bound established above.

Inserting \eqref{int-1} and \eqref{int-2} into \eqref{key-inequality-12}, we conclude
\begin{equation*}
\int_0^\infty \Capp[\Omega]{p}(M_{at},M_t^c) d(t^p)  
\leq \left[\frac{2p\ln a}{(a-1)^p}+2\left(a^\frac{p}{p-1}-1\right)^{1-p}\right]\sum\limits_{\{x,y\}\in E(\Omega,\overline{\Omega})}w_{xy} |f(x)-f(y)|^p.
\end{equation*}
Let $C(a,p)=\frac{2p\ln a}{(a-1)^p}+2\left(a^\frac{p}{p-1}-1\right)^{1-p}$. We complete the proof.
\end{proof}

\section{The first Dirichlet eigenvalues of the \texorpdfstring{$p$}{p}-Laplacian}\label{Dirichlet case}

In this section, we prove Theorem~\ref{capacity-est-Dcase} and Theorem~\ref{Dirichlet-eigenvalue-capacity-estimate-on-infinite-graph}, which give two-sided bounds for $\lambda_{1,p}(\Omega)$ in terms of the $p$-isocapacitary constant $\alpha_p^D(\Omega)$. The upper bound is obtained variationally via a near-optimal test function, and the lower bound follows from the coarea formula applied to the first eigenfunction.

\begin{proof}[Proof of Theorem \ref{capacity-est-Dcase}]
We first show that $\lambda_{1,p}(\Omega)\leq \alpha_p^D(\Omega).$ Recall that 
\[
\alpha_{p}^D(\Omega):=\inf\limits_{A\subset \Omega} \frac{\Capp[\Omega]{p}(A,\delta\Omega)}{m(A)}.
\]
Let $A\subset \Omega$ be a finite subset such that 
\[\alpha_{p}^D(\Omega)=\frac{\Capp[\Omega]{p}(A,\delta\Omega)}{m(A)}.\]
Then there exists a function $f$ satisfying $f|_A=1$ and $f|_{\delta\Omega}=0$, such that $\mathcal{E}^\Omega_p(f,f)=\Capp[\Omega]{p}(A,\delta\Omega)$. By characterization of the Rayleigh quotient, we have 
\[ \lambda_{1,p}(\Omega)\leq \frac{\mathcal{E}_p^\Omega (f,f)}{\|f\|^p_{p,\Omega}}.
\]
Together with  
\[
\|f\|_{p,\Omega}^p=\sum_{x\in\Omega}|f(x)|^pm(x)\geq \sum_{x\in A} m(x)=m(A),
\]
yields 
\[
\alpha_{p}^D(\Omega)=\frac{\Capp[\Omega]{p}(A,\delta\Omega)}{m(A)}\geq \frac{\mathcal{E}_p^\Omega (f,f)}{\|f\|_{p,\Omega}^p} \geq \lambda_{1,p}(\Omega). 
\]

To show that $\frac{1}{2^pC(p)}\alpha_p^D(\Omega)\leq \lambda_{1,p}(\Omega)$, let $u>0$ be the first Dirichlet eigenfunction. Then $u= 0$ on $\delta\Omega$. An application of Lemma \ref{coarea formular} and the definition of $\alpha_p^D(\Omega)$ yields 
\begin{align*}
\lambda_{1,p}(\Omega)\|u\|_{p,\Omega}^p
&=\mathcal{E}_p^\Omega(u,u)\\
&\geq \frac{1}{C_p}\int_0^\infty \Capp[\Omega]{p}(\{|u|\geq 2t\},\{|u|<t\}) d(t^p)\\
&\geq \frac{1}{C_p}\int_0^\infty \Capp[\Omega]{p}(\{|u|\geq 2t\},\delta\Omega) d(t^p)\\
&\geq\frac{1}{C_p}\int_0^\infty \alpha_p^D(\Omega)m\left(\{|u|\geq 2t\}\right)d(t^p)\\
&=\frac{\alpha_p^D(\Omega)}{C_p}\int_0^\infty\sum\limits_{x\in \{|u|\geq 2t\}} m(x) d(t^p)\\
&=\frac{\alpha_p^D(\Omega)}{C_p}\sum_{x\in \{u(x)\neq 0\}}m(x)\int_0^{\frac{|u(x)|}{2}}d(t^p)\\
&=\frac{\alpha_p^D(\Omega)}{2^pC_p}\|u\|_{p,\Omega}^p,
\end{align*}
where we used Tonelli's theorem to interchange the sum and the integral in the penultimate equality.
Hence, we obtain
\[
\lambda_{1,p}(\Omega)\geq \frac{\alpha_p^D(\Omega)}{2^pC_p}.
\]
This completes the proof.
\end{proof}

For an infinite weighted graph $G=(V,E,w,m),$ we define the bottom of the spectrum of the $p$-Laplacian on $G$ as 
\begin{equation}
\lambda_{1,p}(G)=\inf\limits_{f\neq 0}\frac{\mathcal{E}_p(f,f)}{\|f\|^p_{l^p,V}}.
\end{equation}
Let $\mathcal{W}=\{W_i\}_{i=1}^\infty$ be an exhaustion of $G$ and $\lambda_{1,p}(W_i)$ denote the first Dirichlet eigenvalue of $p$-Laplacian. By the Rayleigh quotient characterization, for any $i\in \mathbb{N}_+$, 
\[
\lambda_{1,p}(W_i)\geq \lambda_{1,p}(W_{i+1}).
\]
Thus, we conclude that
\begin{equation}
    \lambda_{1,p}(G)=\lim\limits_{i\to \infty}\lambda_{1,p}(W_i).
\end{equation}
For any subsets $A\subset V$ with $|A|<+\infty$, there exists $i\in \mathbb{N}_+$ such that $A\subset W_i$. Let
\begin{equation*}
\Capp[W_i]{p}(A)
=\inf\{\mathcal{E}^{W_i}_p(f,f)\ | f|_A=1, f|_{\delta W_i}=0\}.
\end{equation*}
It follows that
\[
\Capp[W_i]{p}(A)\geq \Capp[W_{i+1}]{p}(A).
\]

\begin{proof}[Proof of Theorem \ref{Dirichlet-eigenvalue-capacity-estimate-on-infinite-graph}]
By Theorem \ref{capacity-est-Dcase}, we have
\begin{equation*}
\frac{1}{2^pC_p}\lim\limits_{i\to\infty}\alpha_p^D(W_i)\leq \lambda_{1,p}(G)\leq \lim\limits_{i\to \infty}\alpha_p^D(W_i).    
\end{equation*}
It is sufficient to show that $\alpha_p^D(G)=\lim\limits_{i\to \infty}\alpha_p^D(W_i)$.
For any $A\subset V$ with $|A|<+\infty$, there exist $i\in\mathbb{N}_+$ such that $A\subset W_i$. 
\begin{equation*}
\frac{\Capp[V]{p}(A)}{m(A)}=\lim\limits_{i\to \infty}\frac{\Capp[W_i]{p}(A)}{m(A)}\geq \lim\limits_{i\to\infty}\sup\alpha_p^D(W_i).  
\end{equation*}
Hence,
\[
\alpha_p^D(G)\geq \lim\limits_{i\to\infty}\sup\alpha_p^D(W_i).
\]

On the other hand, for any $i\in\mathbb{N}_+$, let $A_i\subset W_i$ be a finite subset such that $\alpha_p^D(W_i)=\frac{\Capp[W_i]{p}(A_i)}{m(A_i)}$.
Then $|A_i|<+\infty$ and
\begin{equation*}
\alpha_p^D(W_i)=\frac{\Capp[W_i]{p}(A_i)}{m(A_i)}\geq \lim\limits_{j\to \infty}\frac{\Capp[W_j]{p}(A_i)}{m(A_i)}
=\frac{\Capp[V]{p}(A_i)}{m(A_i)}\geq \alpha_p^D(G).
\end{equation*}
This completes the proof.
\end{proof}

\section{The first nonzero Neumann eigenvalue of the \texorpdfstring{$p$}{p}-Laplacian}\label{The first Neumann eigenvalues of the $p$-Laplacian}
We now consider the Neumann case. While Dirichlet eigenvalues are determined by the relative $p$-capacity between subsets and the boundary, the Neumann eigenvalue $\mu_{1,p}(\Omega)$ depends fundamentally on the internal connectivity and bottleneck structure of the domain. 

\begin{theorem}\label{Neumann-isocap-est-thm5.1}
 Let $G=(V,E,m,w)$ be a weighted graph. Let $\Omega\subset V$ be a finite subset and the cardinality of $\Omega$ is not less than 2. Then 
\begin{equation}\label{Neum-eigen-cap-est}
\frac{1}{2^pC_p}\overline{\alpha}_p^N(\overline{\Omega})\leq \mu_{1,p}(\Omega)\leq 2^{p-1}\overline{\alpha}_p^N(\overline{\Omega}),
\end{equation}
where $C_p$ is given by \eqref{C(p)} and
\[
\overline{\alpha}_p^N(\overline{\Omega})=\inf\limits_{A,B\in\mathcal{P}(\overline{\Omega})}\frac{\Capp[\Omega]{p}(A,B)}{m(A\cap \Omega)\wedge m(B\cap \Omega)}.
\]
\end{theorem}

\begin{proof}
For the Neumann boundary condition, we only consider the subgraph $G_\Omega.$ Without loss of generality, we may assume that $E(\delta\Omega,\delta\Omega)=\emptyset,$ i.e., $w(x,y)=0$ for any $x,y\in \delta\Omega.$

We first show that $\mu_{1,p}(\Omega)\leq 2^{p-1}\overline{\alpha}_p^N(\overline{\Omega})$. Let $A, B\subset \overline{\Omega}$ such that 
\[
\overline{\alpha}_p^N(\overline{\Omega})=\frac{\Capp[\Omega]{p}(A,B)}{m(A\cap\Omega)\wedge m(B\cap\Omega)}.
\]
Then there exists $f$ satisfying $f|_A=1$ and $f|_B=0$ such that
\[
\mathcal{E}_p^\Omega(f,f)=\Capp[\Omega]{p}(A,B)=\overline{\alpha}_p^N(\overline{\Omega})\left[m(A\cap\Omega)\wedge m(B\cap\Omega)\right].
\]
Then for any $c\in\mathbb{R}$, 
\begin{align*}
\|f-c\|_{p,\overline{\Omega}}^p
&\geq \sum\limits_{x\in A} |f(x)-c|^pm(x) +\sum\limits_{x\in B}|f(x)-c|^p m(x)\\
&=|1-c|^p m(A)+|c|^p m(B)\\
&\geq \frac{1}{2^{p-1}}\left(m(A)\wedge m(B)\right)\\
&\geq  \frac{1}{2^{p-1}}\left[m(A\cap\Omega)\wedge m(B\cap\Omega)\right].
\end{align*}
Thus,
\[
\min\limits_{c\in\mathbb{R}}\|f-c\|^p_{p,\overline{\Omega}}\geq \frac{1}{2^{p-1}}\left[m(A\cap\Omega)\wedge m(B\cap\Omega)\right].
\]
Hence, 
\[
\mu_{1,p}(\Omega)\leq \frac{\mathcal{E}_p^\Omega(f,f)}{\min\limits_{c\in\mathbb{R}}\|f-c\|^p_{p,\overline{\Omega}}}
\leq \frac{\Capp[\Omega]{p}(A,B)}{\frac{1}{2^{p-1}}\left[m(A\cap\Omega)\wedge m(B\cap\Omega)\right]}=2^{p-1}\overline{\alpha}_p^N(\overline{\Omega}).
\] 

Next, we show that $\mu_{1,p}(\Omega)\geq \frac{1}{2^pC_p}\overline{\alpha}_p^N(\overline{\Omega})$. Let $u\in\mathbb{R}^{\overline{\Omega}}$ be a first eigenfunction of $\mu_{1,p}(\Omega)$. As the function changes sign, the sets $\{u>0\}$ and $\{u<0\}$ are not empty. Without loss of generality, we assume 
\[m\left(\{u>0\}\right)\leq \frac{1}{2}m(\overline{\Omega}).\]  
Set $u_+:=u\vee 0$. Let \[\overline{\Omega}_+=\{x\in\overline{\Omega}: u(x)\geq 0\},\ \ \overline{\Omega}_-=\{x\in\overline{\Omega}: u(x)<0\}. \]
Then $u(y)=u_+(y)$ for $y\in \overline{\Omega}_+$ and 
\[
|u(y)-u_+(x)|\geq |u_+(y)-u_+(x)|, \ \ u_+(x)-u(y)\geq u_+(x)-u_+(y)\geq0
\]
for $y\in \overline{\Omega}_-$. According to the Laplacian equation \eqref{Neumann-Laplacian-equation} and boundary condition \eqref{Neumann-bdy-cond}, we have 
\begin{align*} 
\mu_{1,p}(\Omega)\|u_+\|_{p,\Omega}^p
&=\mu_{1,p}(\Omega)\langle |u|^{p-2}u, u_+\rangle_{\Omega}\\
&=\langle -\Delta_p^{G_\Omega} u, u_+\rangle_{\Omega}=\langle -\Delta_p^{G_\Omega} u, u_+\rangle_{{\Omega}}+\left\langle |\nabla u|^{p-2}\frac{\partial u}{\partial n}, u_+\right\rangle_{\delta\Omega}\\
&=\sum\limits_{x\in \overline{\Omega}}\left(\sum\limits_{y\in\overline{\Omega}_+\cup\overline{\Omega}_-}w_{xy}|u(y)-u_+(x)|^{p-2}\left(u_+(x)-u(y)\right)\right)u_+(x)\\
&\geq\sum\limits_{x\in\overline\Omega}\left(\sum\limits_{y\in\overline{\Omega}}w_{xy}|u_+(y)-u_+(x)|^{p-2}\left(u_+(x)-u_+(y)\right)\right)u_+(x)\\
&=\mathcal{E}_p^\Omega(u_+,u_+).
\end{align*}

By Lemma \ref{coarea formular}, we have 
\begin{align*}
\mathcal{E}_p^\Omega(u_+,u_+) 
&\geq \frac{1}{C_p}\int_0^\infty  \Capp[\Omega]{p}\left(\{u_+\geq 2t\}, \{u_+< t\}\right)pt^{p-1}dt\\
&\geq \frac{1}{C_p}\int_0^\infty  \Capp[\Omega]{p}\left(\{u_+\geq 2t\}, \{u_+\leq 0\}\right)pt^{p-1}dt\\
&\geq  \frac{\overline{\alpha}_p^N(\overline{\Omega})}{C_p}\int_0^\infty m\left(\{u_+\geq 2t\}\cap\Omega\right)\wedge m\left(\{u_+\leq 0\}\cap\Omega\right) pt^{p-1}dt\\
&=\frac{\overline{\alpha}_p^N(\overline{\Omega})}{C_p }\int_0^\infty m\left(\{u_+\geq 2t\}\cap\Omega\right)pt^{p-1}dt\\
&=\frac{\overline{\alpha}_p^N(\overline{\Omega})}{C(p)} \int_0^\infty \sum\limits_{x\in \{u_+(x)\geq 2t\}\cap\Omega} m(x)pt^{p-1}dt\\
&=\frac{\overline{\alpha}_p^N(\overline{\Omega})}{C_p}\sum\limits_{x\in \{|u_+|\geq 0\} \cap\Omega } \int_0^{\frac{|u_+(x)|}{2}} pt^{p-1}dt
=\frac{\overline{\alpha}_p^N(\overline{\Omega})}{2^pC_p}\|u_+\|_{p,{\Omega}}^p,
\end{align*}
where the third equality follows again from Tonelli's theorem. Hence,
\[
\mu_{1,p}(\Omega)\geq \frac{\overline{\alpha}_p^N(\overline{\Omega})}{2^pC_p}.
\]
\end{proof}

\begin{proof}[Proof of Theorem \ref{Neumann-capacity-estimate}] 
By Theorem \ref{Neumann-isocap-est-thm5.1}, it is sufficient to prove $\alpha_p^N(\Omega)=\overline{\alpha}_p^N(\overline{\Omega})$.
On the one hand, for any $\{A, B\}\in\mathcal{P}(\overline{\Omega})$, we have 
\begin{equation*}
\frac{\Capp[\Omega]{p}(A,B)}{m(A\cap\Omega)\wedge m(B\cap \Omega)}\geq \frac{\Capp[\Omega]{p}(A\cap\Omega,B\cap\Omega)}{m(A\cap\Omega)\wedge m(B\cap \Omega)}\geq \alpha_p^N(\Omega).  
\end{equation*}
Then $\overline{\alpha}_p^N(\overline{\Omega})\geq \alpha_p^N(\Omega)$. 
On the other hand, let $\{A, B\}\in\mathcal{P}(\Omega)$ such that 
\begin{align*}
\alpha_p^N(\Omega)
&=\frac{\Capp[\Omega]{p}(A,B)}{m(A)\wedge m(B)}\\
&=\frac{\Capp[\Omega]{p}(A,B)}{m(A\cap\Omega)\wedge m(B\cap\Omega)}\\
&\geq \inf\limits_{A,B\in\mathcal{P}(\overline{\Omega})}\frac{\Capp[\Omega]{p}(A,B)}{m(A\cap\Omega)\wedge m(B\cap\Omega)}\\
&=\overline{\alpha}_p^N(\overline{\Omega}).
\end{align*}
This concludes the proof.
\end{proof}

\section{The first nonzero Steklov eigenvalue of the \texorpdfstring{$p$}{p}-Laplacian}\label{First Steklov Eigenvalue of $p$-Laplacian}

In this section, we establish two-sided estimates for the first nonzero eigenvalue using the notion of $p$-isocapacitary. The key idea is to construct a sequence of measures on $G^\Omega=(\overline{\Omega}, E, w,m)$ such that the corresponding first nonzero eigenvalues and $p$-isocapacitary converge to the first nonzero Steklov eigenvalue for the $p$-Laplacian and the Steklov type isocapacitary, respectively.

Let $G=(V,E,w,m)$ be a weighted graph and let $\Omega\subset V$ be a finite subset with the boundary of the vertices such that $|\delta\Omega|\geq 2$. Define $\sigma_{1,p}(\Omega)$ as the first nonzero $p$-Steklov eigenvalue on $\Omega$. According to the Rayleigh quotient characterization, we have
\begin{equation}\label{steklov-eigen-definition}
 \sigma_{1,p}(\Omega)=\inf\limits_{\substack{0\neq u\in \mathbb{R}^{\delta\Omega}} }\frac{\mathcal{E}^{\Omega}_p(u_h, u_h)}{\min\limits_{c\in\mathbb{R}}\|u-c\|_{p,\delta\Omega}^p},   
\end{equation}
where $u_h$ is the unique $p$-harmonic extension of $u$ to $\Omega$. 
It is the unique minimizer of the strictly convex energy 
$\mathcal{E}_p^\Omega(\cdot,\cdot)$ on the affine set of functions 
agreeing with $u$ on $\delta\Omega$, which is a nonempty closed subset 
of $\mathbb{R}^{\overline{\Omega}}$.

Our main results are as follows. 
\begin{theorem}\label{Steklov-isocap-est-thm6.1}
 Let $G=(V,E,w,m)$ be a weighted graph. Let $\Omega\subset V$ be a finite subset with the cardinality $|\delta\Omega|$ is not less than $2$. Then 
\begin{equation}\label{Steklov-eigen-cap-est}
\frac{1}{2^pC_p}\overline{\alpha}_p^S(\overline{\Omega})\leq \sigma_{1,p}(\Omega)\leq 2^{p-1}\overline{\alpha}_p^S(\overline{\Omega}),
\end{equation}
where $C_p$ is given by \eqref{C(p)} and
\[
\overline{\alpha}_p^S(\overline{\Omega})=\inf\limits_{A,B\in \mathcal{P}(\overline{\Omega})}\frac{\Capp[\Omega]{p}(A,B)}{m(A\cap \delta\Omega)\wedge m(B\cap \delta\Omega)}.
\]
\end{theorem}

Next, we construct a sequence of new finite graphs $\{G^{(k)}\}_{k=1}^\infty$, where $G^{(k)}=(\overline{\Omega}, E(\Omega,\overline{\Omega}),w,\\ m^{(k)})$. The measure $m^{(k)}$ is defined such that it satisfies 
\begin{align}
\begin{cases}\label{weighted-S-cases}
m^{(k)}|_{\delta\Omega}=m;\\
m^{(k)}|_\Omega=\frac{m}{k}. 
\end{cases}
\end{align}

We now prove the convergence of the first nonzero
eigenvalue of the $p$-Laplacian under the sequence of measures. 
\begin{lemma}\label{p-Steklov-p-Neumann-Laplacian}
Let $\mu_{1,p}(G^{(k)})$ be the first nonzero eigenvalue of the $p$-Laplacian on $G^{(k)}$ with the measure \eqref{weighted-S-cases}. Then 
\[
\lim\limits_{k \to \infty} \mu_{1,p}(G^{(k)})= \sigma_{1,p}(\Omega).
\]
\end{lemma}

\begin{proof}
On the one hand, let $\phi\in\mathbb{R}^{\delta\Omega}$ such that
\[
\sigma_{1,p}(\Omega)=\frac{\mathcal{E}_p^\Omega(\Phi,\Phi)}{\min\limits_{c\in\mathbb{R}}\|\phi-c\|_{p,\delta\Omega}^p},
\]
where $\Phi$ is the $p$-harmonic extension of $\phi$.  Define constant $c_k$ via
\[
\sum\limits_{x\in\overline{\Omega}} |\Phi(x)- c_k|^{p-2}(\Phi(x)- c_k)m^{(k)}(x)=0.
\]
Then 
\begin{align*}
\mu_{1,p}(G^{(k)}) 
&\leq \frac{\mathcal{E}_p^\Omega(\Phi -c_k, \Phi-c_k)}{\sum\limits_{x\in\overline{\Omega}}|\Phi(x)-c_k|^pm^{(k)}(x)} \\
&=\frac{\mathcal{E}_p^\Omega(\Phi, \Phi)}{\sum\limits_{x\in\overline{\Omega}}|\Phi(x)-c_k|^pm^{(k)}(x)} \\
&=\frac{\sigma_{1,p}(\Omega)\min\limits_{c\in\mathbb{R}}\|\phi-c\|^p_{p,\delta\Omega}}{\frac{1}{k}\sum\limits_{x\in\Omega}|\Phi(x)-c_k|^pm(x)+\sum\limits_{x\in\delta\Omega}|\phi(x)-c_k|^pm(x)}\\
&\leq\sigma_{1,p}(\Omega).
\end{align*}
Hence,
\[
 \limsup\limits_{k\to \infty}  \mu_{1,p}(G^{(k)}) \leq \sigma_{1,p}(\Omega).
\]

On the other hand, let $u_k$ be an eigenfunction corresponding to  $\mu_{1,p}(G^{(k)})$ satisfying
\[
\mu_{1,p}(G^{(k)}) = \frac{\mathcal{E}_p^\Omega(u_k,u_k)}{\sum\limits_{x\in\overline{\Omega}}|u_k(x)|^pm^{(k)}(x)} 
\]
and
\[
\sum_{x\in\overline{\Omega}} |u_k(x)|^{p-2} u_k(x)  m^{(k)}(x) = 0.
\]
We additionally assume that
\[
\| u_k \|_{l^{\infty}, \overline{\Omega}} = 1.
\]
Note that
\[
\limsup\limits_{k\to \infty}  \mu_{1,p}(G^{(k)})<+\infty. 
\] 
For each fixed $k$, the function space is isomorphic to 
$\mathbb{R}^{|G^{(k)}|}$ and the given sequences are bounded. 
By the Bolzano--Weierstrass theorem, we can extract a convergent 
subsequence for each $k$. A diagonal argument then yields an 
increasing sequence of positive integers $\{k_n\}$ with 
$k_n\to\infty$, a nonnegative number $\mu_{1,p}^{(\infty)}$, 
and a function $u_\infty$ such that
\[
\lim_{n\to\infty} \mu_{1,p}(G^{(k_n)}) = \mu_{1,p}^{(\infty)},
\qquad
\lim_{n\to\infty} \|u_{k_n} - u_\infty\|_{l^{\infty},\overline{\Omega}} = 0.
\]
Since the $l^{\infty}$-norm is continuous under $l^{\infty}$-convergence, 
the limit function inherits the normalization, i.e.\ 
$\|u_\infty\|_{l^{\infty},\overline{\Omega}} = 1$.

Recall that each $G^{(k_n)}$ is a graph without boundary whose vertex 
set is $\overline{\Omega}$, hence $\Delta_p^{G^{(k_n)}}$ is defined on 
all of $\overline{\Omega}$ and the eigenvalue equation holds pointwise 
on $\overline{\Omega}$.
All function spaces appearing here are finite-dimensional, so  convergence in $l^{\infty}$ implies convergence in $l^p$:
\[
\lim_{n\to\infty}\|u_{k_n}-u_\infty\|_{l^p,\overline{\Omega}}=0.
\]
In finite dimensions, $l^{\infty}$-convergence implies pointwise convergence as well, so
$u_{k_n}(x)\to u_\infty(x)$ for every $x\in\overline{\Omega}$.
By continuity of $t\mapsto|t|^{p-2}t$, this yields
$|u_{k_n}|^{p-2}u_{k_n}\to|u_\infty|^{p-2}u_\infty$ pointwise on $\overline{\Omega}$.

We now pass to the limit $n\to\infty$ in the relations satisfied by 
$u_{k_n}$.

\begin{itemize}
\item[(i)] For $x\in\delta\Omega$, the identity
  \[
  |\nabla u_{k_n}|^{p-2}\frac{\partial u_{k_n}}{\partial n}(x)
  =-\Delta_p^{G^{(k_n)}}u_{k_n}(x)
  =\mu_{1,p}(G^{(k_n)})|u_{k_n}(x)|^{p-2}u_{k_n}(x)
  \]
  together with the pointwise convergence of $u_{k_n}$ and 
  $\mu_{1,p}(G^{(k_n)})\to\mu_{1,p}^{(\infty)}$ yields
  \[
  |\nabla u_\infty|^{p-2}\frac{\partial u_\infty}{\partial n}(x)
  =\mu_{1,p}^{(\infty)}|u_\infty(x)|^{p-2}u_\infty(x).
  \]

\item[(ii)] For $x\in\Omega$, we have
  \[
  -k_n\Delta_p^{G^{(1)}}u_{k_n}(x)
  =\mu_{1,p}(G^{(k_n)})|u_{k_n}(x)|^{p-2}u_{k_n}(x).
  \]
  Dividing by $k_n$ and passing to the limit yields
  \[
  -\Delta_p^{G^{(1)}}u_\infty(x)=0.
  \]

\item[(iii)] The orthogonality condition for a non-constant eigenfunction gives
  \[
  \sum_{x\in\overline{\Omega}}|u_{k_n}(x)|^{p-2}u_{k_n}(x)\,m^{(k_n)}(x)=0.
  \]
  In the limit $n\to\infty$ we obtain
  \[
  \sum_{x\in\delta\Omega}|u_\infty(x)|^{p-2}u_\infty(x)\,m(x)=0.
  \]
\end{itemize}
Thus, $u_\infty$ is a Steklov eigenfunction corresponding to 
$\mu_{1,p}^{(\infty)}$, and in particular $\mu_{1,p}^{(\infty)}\neq0$. 
Taking \[\mu_{1,p}^{(\infty)}=\liminf\limits_{k\to\infty}\mu_{1,p}(G^{(k)}),\] 
we conclude
\[
\sigma_{1,p}(\Omega)\leq\mu_{1,p}^{(\infty)}
=\liminf_{k\to\infty}\mu_{1,p}(G^{(k)}).
\]

Therefore, we obtain
\[
\lim\limits_{k \to \infty} \mu_{1,p}(G^{(k)})=\sigma_{1,p}(\Omega).
\]
This completes the proof.
\end{proof}

To prove Theorem \ref{Steklov-isocap-est-thm6.1}, we need to show that $\lim\limits_{k\to\infty}\alpha_p(G^{(k)})$ is well-defined.

\begin{lemma}\label{mono-bdd}
For any $\Omega\subset V$ with $|\delta\Omega|\geq 2$ and any $i\in\mathbb{N}_+$, we have
\[
0<\alpha_p(G^{(i)})\leq \alpha_p(G^{(i+1)})\leq C,
\]
where $C$ is a constant depending only on $G_\Omega$.
\end{lemma}

\begin{proof}
Let $A,B\in\mathcal{P}(\overline{\Omega})$ be finite subsets such that 
\[
\alpha_p(G^{(i+1)})=\frac{\Capp[\Omega]{p}(A,B)}{m^{(i+1)}(A)\wedge m^{(i+1)}(B)}.
\]
Since $m^{(i+1)}\leq m^{(i)}$, it follows that
\[
m^{(i+1)}(A)\wedge m^{(i+1)}(B)\leq m^{(i)}(A)\wedge m^{(i)}(B).
\]
This implies that
\[
\alpha_p(G^{(i+1)})\geq \frac{\Capp[\Omega]{p}(A,B)}{ m^{(i)}(A)\wedge m^{(i)}(B)}\geq \alpha_p(G^{(i)}).
\]
Similarly, we have $\alpha_p(G^{(i)})\geq \alpha_p(\overline{\Omega})>0$.

Moreover, since $|\delta\Omega|\geq 2$, we can choose $x_1,x_2\in\delta\Omega$. For any $i\in\mathbb{N}_+$,  it follows that 
\begin{equation*}
\alpha_p(G^{(i)})\leq \frac{\Capp[\Omega]{p}(\{x_1\}, \{x_2\})}{m^{(i)}(x_1)\wedge m^{(i)}(x_2)}=\frac{\Capp[\Omega]{p}(\{x_1\}, \{x_2\})}{m(x_1)\wedge m(x_2)}:=C.    
\end{equation*}
The proof is complete.
\end{proof}

\begin{proof}[Proof of Theorem \ref{Steklov-isocap-est-thm6.1}]
Combining Lemma \ref{p-Steklov-p-Neumann-Laplacian} and Lemma \ref{mono-bdd}, we obtain
\[
\frac{1}{2^pC_p}\lim\limits_{k\to\infty}\alpha_p(G^{(k)})\leq \sigma_{1,p}(\Omega)\leq 2^{p-1}\lim\limits_{k\to\infty}\alpha_p(G^{(k)})   
\]
with $C_p$ is given by \eqref{C(p)}.
Thus, it is sufficient to prove that $\overline{\alpha}_p^S(\overline{\Omega})=\lim\limits_{k\to\infty}\alpha_p(G^{(k)})$.

Recall that  
\[
\overline{\alpha}_p^S(\overline{\Omega})=\inf\limits_{A,B\subset \overline{\Omega}}\frac{\Capp[\Omega]{p}(A,B)}{m(A\cap\delta\Omega)\wedge m(B\cap\delta\Omega)}.
\] 
For any $A,B\subset\overline{\Omega}$ and $k\in\mathbb{N}_+$, we have 
\[m(A\cap \delta\Omega)=m^{(k)}(A\cap \delta\Omega)\leq m^{(k)}(A), m(B\cap \delta\Omega)=m^{(k)}(B\cap \delta\Omega)\leq m^{(k)}(B),\]
thus,
\[
\frac{\Capp[\Omega]{p}(A,B)}{m(A\cap \delta\Omega)\wedge m(B\cap \delta\Omega)}
\geq\frac{\Capp[\Omega]{p}(A,B)}{m^{(k)}(A)\wedge m^{(k)}(B)}
\geq \alpha_p(G^{(k)}).
\]
This implies that
\[
\overline{\alpha}_p^S(\overline{\Omega})\geq \lim\limits_{k\to\infty}\alpha_p(G^{(k)}).
\]

On the other hand, for any $k\in \mathbb{N}_+$, let $A^{(k)},B^{(k)}\subset \overline{\Omega}$ be finite subsets such that 
\[
\alpha_p(G^{(k)})=\frac{\Capp[\Omega]{p}(A^{(k)},B^{(k)})}{m^{(k)}(A^{(k)})\wedge m^{(k)}(B^{(k)})}.
\]
Since $\overline{\Omega}$ is finite, we can choose a subsequence of $\{A^{(k)}\},\{B^{(k)}\}$, still denoted by $A^{(k)}$ and $B^{(k)}$, such that 
\[
A^{(k)}=A\subset\overline{\Omega}, B^{(k)}=B\subset\overline{\Omega}.
\]
Thus, we have 
\begin{align*}
\lim\limits_{k\to \infty}\alpha_p(G^{(k)})
&=\lim\limits_{k\to \infty}\frac{\Capp[\Omega]{p}(A,B)}{m^{(k)}(A)\wedge m^{(k)}(B)}\\
&=\frac{\Capp[\Omega]{p}(A, B)}{m^{(k)}(A\cap\delta\Omega)\wedge m^{(k)}(B\cap\delta\Omega)}\\
&\geq \overline{\alpha}^S_p(\overline{\Omega}).
\end{align*}
Hence, we obtain $\overline{\alpha}_p^S(\overline{\Omega})=\lim\limits_{k\to\infty}\alpha_p(G^{(k)})$.

\end{proof}

\begin{proof}[Proof of Theorem \ref{Steklov-capacity-estimate}] 
By Theorem \ref{Steklov-isocap-est-thm6.1}, it is sufficient to prove $\alpha_p^S(\Omega)=\overline{\alpha}_p^S(\overline{\Omega}).$

On one hand, for any $A, B\subset \overline{\Omega}$, we have 
\begin{equation*}
\frac{\Capp[\Omega]{p}(A,B)}{m(A\cap\delta\Omega)\wedge m(B\cap \delta\Omega)}\geq \frac{\Capp[\Omega]{p}(A\cap\delta\Omega,B\cap \delta\Omega)}{m(A\cap\delta\Omega)\wedge m(B\cap \delta\Omega)}\geq \alpha_p^S(\Omega),    
\end{equation*}
then $\overline{\alpha}_p^S(\overline{\Omega})\geq \alpha_p^N(\Omega)$. 
On the other hand, let $A, B\subset \delta\Omega$ such that 
\begin{align*}
\alpha_p^S(\Omega)
&=\frac{\Capp[\Omega]{p}(A,B)}{m(A)\wedge m(B)}\\
&=\frac{\Capp[\Omega]{p}(A,B)}{m(A\cap\delta\Omega)\wedge m(B\cap\delta\Omega)}\\
&\geq \inf\limits_{A,B\subset \overline{\Omega}}\frac{\Capp[\Omega]{p}(A,B)}{m(A\cap\delta\Omega)\wedge m(B\cap\delta\Omega)}\\
&=\overline{\alpha}_p^S(\overline{\Omega}).
\end{align*}
This concludes the proof.
\end{proof}

We present a direct proof of Theorem~\ref{Steklov-capacity-estimate} based on the coarea formula (Lemma~\ref{coarea formular}). 
\begin{proof}[Alternative proof of Theorem \ref{Steklov-capacity-estimate}]
We first prove the upper bound. Let $\{A,B\}\in\mathcal{P}(\delta\Omega)$ and let $f$ be a minimizer for $\Capp[\Omega]{p}(A,B)$. Then
\[
f|_A=1,\qquad f|_B=0,\qquad
\mathcal{E}_p^\Omega(f,f)=\Capp[\Omega]{p}(A,B).
\]
Since $f$ is prescribed only on $A\cup B\subset\delta\Omega$, the Euler--Lagrange equation gives
\[
\Delta_p f=0\quad\text{in }\Omega.
\]
Thus $f$ is the $p$-harmonic extension of $f|_{\delta\Omega}$. Set $\varphi=f|_{\delta\Omega}$. For every $c\in\mathbb R$,
\[
\begin{aligned}
\|\varphi-c\|_{p,\delta\Omega}^p
&\geq |1-c|^p m(A)+|c|^p m(B)\\
&\geq \left(m(A)\wedge m(B)\right)\left(|1-c|^p+|c|^p\right)\\
&\geq 2^{1-p}\left(m(A)\wedge m(B)\right).
\end{aligned}
\]
By the variational characterization \eqref{steklov-eigen-definition},
\[
\sigma_{1,p}(\Omega)
\leq \frac{\mathcal{E}_p^\Omega(f,f)}
{\min\limits_{c\in\mathbb R}\|\varphi-c\|_{p,\delta\Omega}^p}
\leq 2^{p-1}\frac{\Capp[\Omega]{p}(A,B)}{m(A)\wedge m(B)}.
\]
Taking the infimum over all $\{A,B\}\in\mathcal{P}(\delta\Omega)$ gives
\[
\sigma_{1,p}(\Omega)\leq 2^{p-1}\alpha_p^S(\Omega).
\]

For the lower bound, let $u$ be an eigenfunction corresponding to the first nonzero Steklov eigenvalue $\sigma_{1,p}(\Omega)$. Applying Lemma~\ref{lemma:green-formula} with $f=u$ and $g=1$, we get
\[
0=\mathcal{E}_p^\Omega(u,1)=\left\langle |\nabla u|^{p-2}\frac{\partial u}{\partial n},1\right\rangle_{\delta\Omega}=\sigma_{1,p}(\Omega)\sum_{x\in\delta\Omega}|u(x)|^{p-2}u(x)m(x).
\]
Thus,
\[
\sum_{x\in\delta\Omega}|u(x)|^{p-2}u(x)m(x)=0.
\]
Which implies
\begin{equation}\label{eq:boundary-sign-change}
m(\{u>0\}\cap\delta\Omega)>0
\quad\text{and}\quad
m(\{u<0\}\cap\delta\Omega)>0.
\end{equation}
Replacing $u$ by $-u$ if necessary, we may assume that
\begin{equation}\label{eq:positive-measure-bound}
m(\{u>0\}\cap\delta\Omega)\leq \frac{1}{2}m(\delta\Omega).
\end{equation}
Set $u^+=u\vee0$. Applying Lemma~\ref{lemma:green-formula} with $f=u$ and $g=u^+$ yields
\begin{equation}\label{eq:steklov-energy-positive}
\sigma_{1,p}(\Omega)\|u^+\|_{p,\delta\Omega}^p
=\mathcal{E}_p^\Omega(u,u^+).
\end{equation}
For $s\in\mathbb R$, write $s^+=s\vee 0$. For every $a,b\in\mathbb R$, since $(a-b)(a^+-b^+)\geq0$ and $|a^+-b^+|\leq |a-b|$, we have
\[
|a-b|^{p-2}(a-b)(a^+-b^+)
\geq |a^+-b^+|^p.
\]
Applying this pointwise inequality to the definition of $\mathcal{E}_p^\Omega$, together with Lemma~\ref{coarea formular}, we have
\begin{equation}\label{eq:energy-coarea}
\mathcal{E}_p^\Omega(u,u^+)
\geq \mathcal{E}_p^\Omega(u^+,u^+)
\geq \frac{1}{C_p}\int_0^\infty\Capp[\Omega]{p}
\left(\{u^+\ge 2t\},\{u^+\le 0\}\right)
\,p t^{p-1}\,dt.
\end{equation}

For $t>0$, set
\[
A_t=\{u^+\ge 2t\}\cap\delta\Omega,\qquad
B=\{u^+\le 0\}\cap\delta\Omega.
\]
Since $A_t\subset \{u>0\}\cap\delta\Omega$, $B=\{u\le 0\}\cap\delta\Omega$ and \eqref{eq:positive-measure-bound},  
\[
m(A_t)\leq \frac{1}{2}m(\delta\Omega),\, \,  m(B)\geq \frac{1}{2}m(\delta\Omega).
\]
By \eqref{eq:boundary-sign-change}, $B\neq \emptyset$. If $A_t\neq\emptyset$,  
\begin{equation}\label{eq:capacity-bound}
\begin{aligned}
\Capp[\Omega]{p}\left(\{u^+\ge 2t\},\{u^+\le 0\}\right)
&\geq \Capp[\Omega]{p}(A_t,B) \\
&\geq \alpha_p^S(\Omega)\left(m(A_t)\wedge m(B)\right)
= \alpha_p^S(\Omega)m(A_t).
\end{aligned}
\end{equation}
If $A_t=\emptyset$, then $m(A_t)=0$, and the desired inequality is trivial. Combining \eqref{eq:steklov-energy-positive}, \eqref{eq:energy-coarea} and \eqref{eq:capacity-bound}, we have
\[
\begin{aligned}
\sigma_{1,p}(\Omega)\|u^+\|_{p,\delta\Omega}^p
&\geq \frac{\alpha_p^S(\Omega)}{C_p}\int_0^\infty m(A_t)\,p t^{p-1}\,dt\\ 
&=\frac{\alpha_p^S(\Omega)}{C_p}
\int_0^\infty m\left(\{u^+\geq 2t\}\cap\delta\Omega\right)\,p t^{p-1}\,dt \\
&=\frac{\alpha_p^S(\Omega)}{2^p C_p} \|u^+\|_{p,\delta\Omega}^p,
\end{aligned}
\]
where we used Tonelli's theorem in the last equality. By \eqref{eq:boundary-sign-change}, $\|u^+\|_{p,\delta\Omega}^p\neq 0$. Therefore,
\[
\sigma_{1,p}(\Omega) \geq \frac{\alpha_p^S(\Omega)}{2^pC_p}.
\]
This completes the proof.
\end{proof}

We present an example to demonstrate the sharpness of the upper bound in the theorem  \ref{Steklov-capacity-estimate}.
\begin{example}[A path graph]
Let $G = (\mathbb{Z}, E)$ be the path graph with unit edge weights and unit vertex weights, $\Omega = \{1, 2, \dots, n-1\}$ for $n \geq 2$. For $p>1$, we have \[
\sigma_{1,p}(\Omega)=2^{p-1}\alpha_p^S(\Omega).
\]
\end{example}

\begin{proof}
Since $\delta\Omega=\{0,n\}$, by \eqref 
{def-alpha-S-isocap}, we can choose $A = \{0\}, B = \{n\}$.  Thus, $m(A) \wedge m(B) = 1$, and
\[
\Capp[\Omega]{p}(A,B)=\inf\left\{\sum\limits_{i=1}^n|f(i)-f(i-1)|^p: f(0)=1,\  f(n)=0,\  f\in l_0(\overline{\Omega})\right\}.
\]
The infimum is attained by $f: \overline{\Omega} \to \mathbb{R}$ satisfying 
\[
f(x) = 1 - x/n\]
for $x \in \{0, 1, \dots, n\}$.  Hence,
    \[
    \alpha_p^S(\Omega) = \inf_{A',B' \subset \delta\Omega} \frac{\operatorname{Cap}_p^\Omega(A',B')}{m(A') \wedge m(B')} = n^{1-p}.
    \]

Next, we consider $\sigma_{1,p}(\Omega)$. By \eqref{steklov-eigen-definition}, 
\[ 
 \sigma_{1,p}(\Omega)=\inf\limits_{0\neq u\in \mathbb{R}^{\delta\Omega}}\frac{\sum\limits_{i=1}^n\left|u_h(i)-u_h(i-1)\right|^p}{\min\limits_{c\in\mathbb{R}}\left(|u(0)-c|^p+|u(n)-c|^p \right)},   
\]
where $u_h$ is the $p$-harmonic extension of $u$ to $\Omega$.

Since $p>1$, $f(c):=|u(0)-c|^p+|u(n)-c|^p$ is a convex function of $c$ and the minimum is attained at $c=\frac{u(0)+u(n)}{2}$. Thus, 
\[
\min\limits_{c\in\mathbb{R}}\left(|u(0)-c|^p+|u(n)-c|^p \right)=2^{1-p}|u(n)-u(0)|^p.
\]

On a path graph, the $p$-harmonic extension is a linear function, explicitly given by $u_h(x) = a + bx$ and $a,b$ are constants. Clearly, 
$u(0)=a$, $u(n)=a+bn$ and $u_h(i)-u_h(i-1)=b$. Then 
\[
\sigma_{1,p}(\Omega) = 2^{p-1} n^{1-p}.
\]
Therefore, $\sigma_{1,p}(\Omega)=2^{p-1}\alpha_p^S(\Omega)$.
\end{proof}

\section{Acknowledgments}
We are grateful to the anonymous referee for suggesting a direct proof of Theorem~\ref{Steklov-capacity-estimate} using the coarea formula. The authors thank Professor Xueping Huang for enlightening discussions and valuable suggestions regarding $p$-capacity. B. H. thanks the University of Jena for its hospitality while this paper was written during a research stay at the University of Jena as part of the Jena Excellence Fellowship Programme. L. W. thanks Tao Wang for helpful discussions and suggestions.  B. H. is supported by NSFC, no.12371056. L. W. is supported by NSFC, no. 12371052, and the Fujian Alliance of Mathematics, no. 2024SXLMMS01.

\section*{Declarations}

\subsection*{Ethical Approval}
Not applicable.

\subsection*{Funding}
B.~H. is supported by NSFC, no.~12371056. 
L.~W. is supported by NSFC, no.~12371052, and the Fujian Alliance of Mathematics, no.~2024SXLMMS01.

\subsection*{Availability of data and materials}
Not applicable.

\bibliography{ref1}

@article {Colbois2024,
    AUTHOR = {Colbois, Bruno and Girouard, Alexandre and Gordon, Carolyn and
              Sher, David},
     TITLE = {Some recent developments on the {S}teklov eigenvalue problem},
   JOURNAL = {Rev. Mat. Complut.},
  FJOURNAL = {Revista Matem\'atica Complutense},
    VOLUME = {37},
      YEAR = {2024},
    NUMBER = {1},
     PAGES = {1--161},
      ISSN = {1139-1138,1988-2807},
   MRCLASS = {58C40 (35P05 35P15 35P20 53A10 58J50)},
  MRNUMBER = {4695859},
       DOI = {10.1007/s13163-023-00480-3},
       URL = {https://doi.org/10.1007/s13163-023-00480-3},
}

@article{Barletta2025,
  author  = {Giuseppina Barletta},
  title   = {L\ensuremath{^\infty} boundedness of the solutions to the {S}chr\"{o}dinger
             equation on noncompact {R}iemannian manifolds},
  journal = {Journal of Differential Equations},
  volume  = {439},
  pages   = {113413},
  year    = {2025},
  doi     = {10.1016/j.jde.2025.113413}
}

@article{Wang-Zhou2024,
  author  = {Jiaxiang Wang and Bin Zhou},
  title   = {Trace inequalities, isocapacitary inequalities, and regularity of the
             complex {H}essian equations},
  journal = {Science China Mathematics},
  volume  = {67},
  number  = {3},
  pages   = {557--576},
  year    = {2024},
  doi     = {10.1007/s11425-022-2100-1}
}

@book{Folland99,
  author    = {Folland, Gerald B.},
  title     = {Real analysis},
  series    = {Pure and Applied Mathematics: A Wiley-Interscience Series of Texts, Monographs and Tracts},
  edition   = {Second},
  note      = {Modern techniques and their applications},
  publisher = {John Wiley \& Sons, Inc., New York},
  year      = {1999},
  pages     = {xiv+386},
  isbn      = {0-471-31716-0},
  mrclass   = {28-01 (26-01 46-01 54-01)},
  mrnumber  = {1683747}
}

@article {1985lambda1,
    AUTHOR = {Alon, N. and Milman, V. D.},
     TITLE = {{$\lambda_1,$} isoperimetric inequalities for graphs, and
              superconcentrators},
   JOURNAL = {J. Combin. Theory Ser. B},
  FJOURNAL = {Journal of Combinatorial Theory. Series B},
    VOLUME = {38},
      YEAR = {1985},
    NUMBER = {1},
     PAGES = {73--88},
      ISSN = {0095-8956,1096-0902},
   MRCLASS = {05C50},
  MRNUMBER = {782626},
       DOI = {10.1016/0095-8956(85)90092-9},
       URL = {https://doi.org/10.1016/0095-8956(85)90092-9},
}

@article {Barlow2009heat,
    AUTHOR = {Barlow, Martin T. and Grigor'yan, Alexander and
              Kumagai, Takashi},
     TITLE = {Heat kernel upper bounds for jump processes and the first exit
              time},
   JOURNAL = {J. Reine Angew. Math.},
  FJOURNAL = {Journal f\"ur die Reine und Angewandte Mathematik. [Crelle's
              Journal]},
    VOLUME = {626},
      YEAR = {2009},
     PAGES = {135--157},
      ISSN = {0075-4102,1435-5345},
   MRCLASS = {58J65 (35K05 60J75)},
  MRNUMBER = {2492992},
MRREVIEWER = {David\ A.\ Croydon},
       DOI = {10.1515/CRELLE.2009.005},
       URL = {https://doi.org/10.1515/CRELLE.2009.005},
}

@article {Bonder2001,
    AUTHOR = {Bonder, Juli\'an Fern\'andez and Rossi, Julio D.},
     TITLE = {Existence results for the {$p$}-{L}aplacian with nonlinear
              boundary conditions},
   JOURNAL = {J. Math. Anal. Appl.},
  FJOURNAL = {Journal of Mathematical Analysis and Applications},
    VOLUME = {263},
      YEAR = {2001},
    NUMBER = {1},
     PAGES = {195--223},
      ISSN = {0022-247X,1096-0813},
   MRCLASS = {35J60 (35J25 35J65)},
  MRNUMBER = {1864315},
MRREVIEWER = {Alan\ V.\ Lair},
       DOI = {10.1006/jmaa.2001.7609},
       URL = {https://doi.org/10.1006/jmaa.2001.7609},
}

@unpublished{Hua-Munch-Wang,
  AUTHOR = {Hua, Bobo and Münch, Florentin and Wang, Tao},
	TITLE = {Cheeger type inequalities associated with isocapacitary constants on graphs},
   YEAR = {2024},
 note = {arXiv:2406.12583 },
}

@book {Chung1997SpectralGT,
    AUTHOR = {Chung, Fan R. K.},
     TITLE = {Spectral graph theory},
    SERIES = {CBMS Regional Conference Series in Mathematics},
    VOLUME = {92},
 PUBLISHER = {Conference Board of the Mathematical Sciences, Washington, DC;
              by the American Mathematical Society, Providence, RI},
      YEAR = {1997},
     PAGES = {xii+207},
      ISBN = {0-8218-0315-8},
   MRCLASS = {58G99 (05C50 35P05 46N20 47N20)},
  MRNUMBER = {1421568},
MRREVIEWER = {Robert\ Brooks},
}

@article {Dodziuk1984difference,
    AUTHOR = {Dodziuk, Jozef},
     TITLE = {Difference equations, isoperimetric inequality and transience
              of certain random walks},
   JOURNAL = {Trans. Amer. Math. Soc.},
  FJOURNAL = {Transactions of the American Mathematical Society},
    VOLUME = {284},
      YEAR = {1984},
    NUMBER = {2},
     PAGES = {787--794},
      ISSN = {0002-9947,1088-6850},
   MRCLASS = {58G32 (35J05 39A12 53C99)},
  MRNUMBER = {743744},
MRREVIEWER = {J.\ L.\ Kazdan},
       DOI = {10.2307/1999107},
       URL = {https://doi.org/10.2307/1999107},
}

@article {Garcia-Azorero-Manfredi-Peral-2006,
    AUTHOR = {Garcia-Azorero, Jesus and Manfredi, Juan J. and Peral, Ireneo
              and Rossi, Julio D.},
     TITLE = {Steklov eigenvalues for the {$\infty$}-{L}aplacian},
   JOURNAL = {Atti Accad. Naz. Lincei Rend. Lincei Mat. Appl.},
  FJOURNAL = {Atti della Accademia Nazionale dei Lincei. Rendiconti Lincei.
              Matematica e Applicazioni},
    VOLUME = {17},
      YEAR = {2006},
    NUMBER = {3},
     PAGES = {199--210},
      ISSN = {1120-6330,1720-0768},
   MRCLASS = {35J60 (35J65 35J70 47J30)},
  MRNUMBER = {2254067},
MRREVIEWER = {Abdelmajid\ Siai},
       DOI = {10.4171/RLM/463},
       URL = {https://doi.org/10.4171/RLM/463},
}

@book{grigor-isoperimetric,
    AUTHOR = {Grigor'yan, Alexander},
     TITLE = {Isoperimetric inequalities and capacities on {R}iemannian
              manifolds},
 BOOKTITLE = {The {M}az'ya anniversary collection, {V}ol. 1
              ({R}ostock, 1998)},
    SERIES = {Oper. Theory Adv. Appl.},
    VOLUME = {109},
     PAGES = {139--153},
 PUBLISHER = {Birkh\"auser, Basel},
      YEAR = {1999},
      ISBN = {3-7643-6201-4},
   MRCLASS = {31C12 (35P15 58J60)},
  MRNUMBER = {1747869},
MRREVIEWER = {Mario\ Bonk},
       DOI = {10.1007/978-3-0348-8675-8\_9},
       URL = {https://doi.org/10.1007/978-3-0348-8675-8_9},
}

@article {Hassannezhad-miclo20,
    AUTHOR = {Hassannezhad, Asma and Miclo, Laurent},
     TITLE = {Higher order {C}heeger inequalities for {S}teklov eigenvalues},
   JOURNAL = {Ann. Sci. \'Ec. Norm. Sup\'er. (4)},
  FJOURNAL = {Annales Scientifiques de l'\'Ecole Normale Sup\'erieure.
              Quatri\`eme S\'erie},
    VOLUME = {53},
      YEAR = {2020},
    NUMBER = {1},
     PAGES = {43--88},
      ISSN = {0012-9593,1873-2151},
   MRCLASS = {58C40 (58J50)},
  MRNUMBER = {4093440},
MRREVIEWER = {Xiaolong\ Li},
       DOI = {10.24033/asens.2417},
       URL = {https://doi.org/10.24033/asens.2417},
}

@article {Holopainen92,
    AUTHOR = {Holopainen, Ilkka},
     TITLE = {Positive solutions of quasilinear elliptic equations on
              {R}iemannian manifolds},
   JOURNAL = {Proc. London Math. Soc. (3)},
  FJOURNAL = {Proceedings of the London Mathematical Society. Third Series},
    VOLUME = {65},
      YEAR = {1992},
    NUMBER = {3},
     PAGES = {651--672},
      ISSN = {0024-6115,1460-244X},
   MRCLASS = {58G30 (30C65 31C12 35J60 53C20)},
  MRNUMBER = {1182105},
MRREVIEWER = {Hitoshi\ Arai},
       DOI = {10.1112/plms/s3-65.3.651},
       URL = {https://doi.org/10.1112/plms/s3-65.3.651},
}

@article {Hua-Huang2018,
    AUTHOR = {Hua, Bobo and Huang, Yan},
     TITLE = {Neumann {C}heeger constants on graphs},
   JOURNAL = {J. Geom. Anal.},
  FJOURNAL = {Journal of Geometric Analysis},
    VOLUME = {28},
      YEAR = {2018},
    NUMBER = {3},
     PAGES = {2166--2184},
      ISSN = {1050-6926,1559-002X},
   MRCLASS = {05C50 (31C12 39A12)},
  MRNUMBER = {3833789},
MRREVIEWER = {Shiping\ Liu},
       DOI = {10.1007/s12220-017-9899-8},
       URL = {https://doi.org/10.1007/s12220-017-9899-8},
}

@article {HHW2017,
    AUTHOR = {Hua, Bobo and Huang, Yan and Wang, Zuoqin},
     TITLE = {First eigenvalue estimates of {D}irichlet-to-{N}eumann
              operators on graphs},
   JOURNAL = {Calc. Var. Partial Differential Equations},
  FJOURNAL = {Calculus of Variations and Partial Differential Equations},
    VOLUME = {56},
      YEAR = {2017},
    NUMBER = {6},
     PAGES = {Paper No. 178, 21},
      ISSN = {0944-2669,1432-0835},
   MRCLASS = {35R02 (05C90 35J05 35J25 35P15 58J32)},
  MRNUMBER = {3722072},
MRREVIEWER = {Gabriella\ Bretti},
       DOI = {10.1007/s00526-017-1260-3},
       URL = {https://doi.org/10.1007/s00526-017-1260-3},
}

@article {Hua-Wang20,
    AUTHOR = {Hua, Bobo and Wang, Lili},
     TITLE = {Dirichlet {$p$}-{L}aplacian eigenvalues and {C}heeger
              constants on symmetric graphs},
   JOURNAL = {Adv. Math.},
  FJOURNAL = {Advances in Mathematics},
    VOLUME = {364},
      YEAR = {2020},
     PAGES = {106997, 34},
      ISSN = {0001-8708,1090-2082},
   MRCLASS = {05C50},
  MRNUMBER = {4058211},
MRREVIEWER = {Mingqing\ Zhai},
       DOI = {10.1016/j.aim.2020.106997},
       URL = {https://doi.org/10.1016/j.aim.2020.106997},
}

@article {Keller-Mugnolo16,
    AUTHOR = {Keller, Matthias and Mugnolo, Delio},
     TITLE = {General {C}heeger inequalities for {$p$}-{L}aplacians on
              graphs},
   JOURNAL = {Nonlinear Anal.},
  FJOURNAL = {Nonlinear Analysis. Theory, Methods \& Applications. An
              International Multidisciplinary Journal},
    VOLUME = {147},
      YEAR = {2016},
     PAGES = {80--95},
      ISSN = {0362-546X,1873-5215},
   MRCLASS = {05C50 (35J92 35R02)},
  MRNUMBER = {3564720},
       DOI = {10.1016/j.na.2016.07.011},
       URL = {https://doi.org/10.1016/j.na.2016.07.011},
}

@article{Mazya1960,
    AUTHOR = {Maz'ya, V. G.},
     TITLE = {Classes of domains and imbedding theorems for function spaces},
JOURNAL ={Soviet Math. Dokl.},
FJOURNAL ={Soviet Mathematics. Doklady},
volume={1},
year={1960},
     PAGES = {882--885},
ISSN = {0197-6788},
   MRCLASS = {46.38},
  MRNUMBER = {126152},
MRREVIEWER = {J.\ Dieudonn\'e},
}

@book{Mazya1985,
    AUTHOR = {Maz'ya, V. G.},
     TITLE = {\newblock {\em \foreignlanguage{russian}{Prostranstva S. L. Soboleva}} (Russian)[Sobolev spaces]},
 PUBLISHER = {Leningrad. Univ., Leningrad},
      YEAR = {1985},
     PAGES = {416},
   MRCLASS = {46E35 (31C15 47F05)},
  MRNUMBER = {807364},
MRREVIEWER = {H.\ Triebel},
}

@article{Mazya1962,
    AUTHOR = {Maz'ya, V. G.},
     TITLE = {The negative spectrum of the higher-dimensional
              {S}chr\"odinger operator},
   JOURNAL = {Dokl. Akad. Nauk SSSR},
  FJOURNAL = {Doklady Akademii Nauk SSSR},
    VOLUME = {144},
      YEAR = {1962},
     PAGES = {721--722},
      ISSN = {0002-3264},
   MRCLASS = {35.80 (47.65)},
  MRNUMBER = {138880},
MRREVIEWER = {Archibald\ Brown},
}

@article {Mazya1964,
    AUTHOR = {Maz'ya, V. G.},
     TITLE = {On the theory of the higher-dimensional {S}chr\"odinger
              operator},
   JOURNAL = {Izv. Akad. Nauk SSSR Ser. Mat.},
  FJOURNAL = {Izvestiya Akademii Nauk SSSR. Seriya Matematicheskaya},
    VOLUME = {28},
      YEAR = {1964},
     PAGES = {1145--1172},
      ISSN = {0373-2436},
   MRCLASS = {35.78 (47.65)},
  MRNUMBER = {174879},
}

@article {Mazya2005Conductor,
    AUTHOR = {Maz'ya, Vladimir},
     TITLE = {Conductor and capacitary inequalities for functions on
              topological spaces and their applications to {S}obolev-type
              imbeddings},
   JOURNAL = {J. Funct. Anal.},
  FJOURNAL = {Journal of Functional Analysis},
    VOLUME = {224},
      YEAR = {2005},
    NUMBER = {2},
     PAGES = {408--430},
      ISSN = {0022-1236,1096-0783},
   MRCLASS = {31C45 (46E35)},
  MRNUMBER = {2146047},
MRREVIEWER = {Noureddine\ A\"issaoui},
       DOI = {10.1016/j.jfa.2004.09.009},
       URL = {https://doi.org/10.1016/j.jfa.2004.09.009},
}

@incollection {mazya2009,
    AUTHOR = {Maz'ya, Vladimir},
     TITLE = {Integral and isocapacitary inequalities},
 BOOKTITLE = {Linear and complex analysis},
    SERIES = {Amer. Math. Soc. Transl. Ser. 2},
    VOLUME = {226},
     PAGES = {85--107},
 PUBLISHER = {Amer. Math. Soc., Providence, RI},
      YEAR = {2009},
      ISBN = {978-0-8218-4801-2; 0-8218-4801-1},
   MRCLASS = {46E35 (35A23)},
  MRNUMBER = {2500513},
       DOI = {10.1090/trans2/226/09},
       URL = {https://doi.org/10.1090/trans2/226/09},
}

@article {Provenzano2022upper,
    AUTHOR = {Provenzano, Luigi},
     TITLE = {Upper bounds for the {S}teklov eigenvalues of the
              {$p$}-{L}aplacian},
   JOURNAL = {Mathematika},
  FJOURNAL = {Mathematika. A Journal of Pure and Applied Mathematics},
    VOLUME = {68},
      YEAR = {2022},
    NUMBER = {1},
     PAGES = {148--162},
      ISSN = {0025-5793,2041-7942},
   MRCLASS = {35P15 (35P30 58J50)},
  MRNUMBER = {4405973},
       DOI = {10.1112/mtk.12119},
       URL = {https://doi.org/10.1112/mtk.12119},
}

@article {Verma2020,
    AUTHOR = {Verma, Sheela},
     TITLE = {Upper bound for the first nonzero eigenvalue related to the
              {$p$}-{L}aplacian},
   JOURNAL = {Proc. Indian Acad. Sci. Math. Sci.},
  FJOURNAL = {Indian Academy of Sciences. Proceedings. Mathematical
              Sciences},
    VOLUME = {130},
      YEAR = {2020},
    NUMBER = {1},
     PAGES = {Paper No. 21, 11},
      ISSN = {0253-4142,0973-7685},
   MRCLASS = {35J92 (35J66 35P15 58J50)},
  MRNUMBER = {4064275},
       DOI = {10.1007/s12044-019-0529-1},
       URL = {https://doi.org/10.1007/s12044-019-0529-1},
}

@article {Yamasaki1977,
    AUTHOR = {Yamasaki, Maretsugu},
     TITLE = {Parabolic and hyperbolic infinite networks},
   JOURNAL = {Hiroshima Math. J.},
  FJOURNAL = {Hiroshima Mathematical Journal},
    VOLUME = {7},
      YEAR = {1977},
    NUMBER = {1},
     PAGES = {135--146},
      ISSN = {0018-2079,2758-9641},
   MRCLASS = {94A20},
  MRNUMBER = {429377},
MRREVIEWER = {C.\ Wang},
       URL = {http://projecteuclid.org/euclid.hmj/1206135953},
}

@article {Yamasaki1986ideal,
    AUTHOR = {Yamasaki, Maretsugu},
     TITLE = {Ideal boundary limit of discrete {D}irichlet functions},
   JOURNAL = {Hiroshima Math. J.},
  FJOURNAL = {Hiroshima Mathematical Journal},
    VOLUME = {16},
      YEAR = {1986},
    NUMBER = {2},
     PAGES = {353--360},
      ISSN = {0018-2079,2758-9641},
   MRCLASS = {31C20},
  MRNUMBER = {855163},
MRREVIEWER = {Moses\ Glasner},
       URL = {http://projecteuclid.org/euclid.hmj/1206130433},
}

\bibliographystyle{plain}

\end{document}